\pdfoutput=1
\documentclass[11pt]{article}

\usepackage[margin=1.1in]{geometry}
\usepackage[T1]{fontenc}
\usepackage{amsmath,amsthm}
\usepackage{lmodern}                   
\usepackage{amssymb}
\usepackage{booktabs}
\usepackage{url}
\usepackage[expansion=false]{microtype}
\usepackage[colorlinks=true,linkcolor=black,citecolor=blue,urlcolor=blue]{hyperref}
\let\plainparagraph\paragraph
\renewcommand{\paragraph}[1]{\plainparagraph[#1]{\boldmath #1}}
\hypersetup{
  pdftitle={New lower bounds for six small hypergraph Ramsey numbers},
  pdfauthor={Christian Lysenst{\o}en},
}

\newtheorem{theorem}{Theorem}
\newtheorem{corollary}[theorem]{Corollary}
\newtheorem{proposition}[theorem]{Proposition}
\newtheorem{lemma}[theorem]{Lemma}
\theoremstyle{remark}
\newtheorem{remark}[theorem]{Remark}

\newcommand{\Z}{\mathbb{Z}}
\newcommand{\G}{G}
\newcommand{\Bfs}{201}

\newcommand{\Bfseven}{401}
\newcommand{\Bsseven}{2001}
\newcommand{\Bsevenseven}{12\,001}

\newcommand{\Beighteight}{588\,001}

\title{\bfseries New lower bounds for six small\\ hypergraph Ramsey numbers}

\author{Christian Lysenst{\o}en\\[4pt]
\small University of Inland Norway\\[2pt]
\small\texttt{christian@lysenstoen.net}}

\date{November 2026}

\begin{document}
\maketitle

\section*{Abstract}
Let $R(k,l;t)$ denote the least $n$ such that every colouring of the
$t$-subsets of an $n$-set with two colours has $k$ points all of whose
$t$-subsets have the first colour, or $l$ points all of whose $t$-subsets have
the second. We give explicit colourings proving
$R(4,5;3)\ge 38$, $R(4,6;3)\ge 71$, $R(5,5;3)\ge 99$,
$R(5,6;4)\ge 201$, $R(6,6;4)\ge 601$ and $R(6,6;5)\ge 84$. The first three
improve the bounds $35$, $63$ and $88$ listed in Radziszowski's dynamic
survey. The last three improve the listed bounds $67$, $133$ and $72$, and
for the two $4$-uniform numbers also the values $159$ and $317$ of a recent
data deposit. The colourings of $3$-subsets are invariant under groups of
affine maps $x\mapsto ax+b$ of $\Z_n$, for $n=37$, $69$ and $97$, in the last
two cases with one further fixed point. Those of $4$-subsets live on the
projective line over $\mathbb{F}_q$, for $q=199$ and $q=599$, and are
invariant under $\mathrm{PSL}(2,q)$ or $\mathrm{PGL}(2,q)$, so each is given
by a table of $q$ entries. For $5$-subsets we show that, for every prime
$q\ge 5$, invariance under $\mathrm{PGL}(2,q)$ forces six points all of whose
$5$-subsets have the same colour; our colouring is invariant under all affine
maps of $\Z_{83}$. By known constructions the bounds give $R(5,5,5;3)\ge 9605$ and
$R_k(5;3)\ge 98^{2^{k-2}}+1$, and $R(6,7;4)\ge 2001$, $R(7,7;4)\ge 12\,001$
and $R(8,8;4)\ge 588\,001$. The colourings, and programs that verify them,
accompany the paper as ancillary files.

\medskip
\noindent\textbf{Keywords:} hypergraph Ramsey numbers; lower bounds;
prescribed automorphisms; projective line; SAT; local search;
computational verification.\\
\textbf{MSC 2020:} 05D10; 05C55; 05C65.

\section{Introduction}\label{sec:intro}

For integers $t\ge 2$ and $k_1,\dots,k_r>t$, the \emph{hypergraph Ramsey
number} $R(k_1,\dots,k_r;t)$ is the least $n$ such that every colouring of the
$t$-subsets of an $n$-set with the colours $1,\dots,r$ has, for some $i$, a
set of $k_i$ points all of whose $t$-subsets have colour $i$. We write
$R_r(k;t)$ when $k_1=\dots=k_r=k$. With two colours, which we call $0$
and~$1$, a colouring of the $t$-subsets of an $n$-set with no $k$ points
whose $t$-subsets all have colour $0$ and no $l$ points whose $t$-subsets all
have colour $1$ is a \emph{$(k,l;t)$-colouring of order $n$}; it proves
$R(k,l;t)\ge n+1$. We call the $3$-subsets \emph{triples}.

Very little is known about these numbers for $t\ge 3$. The only nontrivial
value known is $R(4,4;3)=13$, computed by McKay and
Radziszowski~\cite{mr1}. Radziszowski's dynamic
survey~\cite[Section~7.1]{survey} lists the lower bounds $R(4,5;3)\ge 35$,
due to Dybizba\'nski~\cite{dyb2}, and $R(4,6;3)\ge 63$ and $R(5,5;3)\ge 88$,
also due to Dybizba\'nski~\cite{dyb3}; these come from explicit colourings.
Among the bounds that the survey obtains from recurrence relations and
related methods \cite{sonyl,liu} are $R(5,6;4)\ge 67$, $R(6,6;4)\ge 133$ and
$R(6,6;5)\ge 72$. Two data deposits of 2026 report further bounds:
$R(4,6;3)\ge 64$~\cite{yang}, and $R(5,5;4)\ge 49$ and
$R(5,6;4)\ge 159$~\cite{ozdemir}. The last two come from colourings of the
projective line that are invariant under $\mathrm{PSL}(2,q)$; from the second
of them and a recurrence, \cite{ozdemir} also states $R(6,6;4)\ge 317$,
without a colouring.

\begin{theorem}\label{thm:main}
$R(4,5;3)\ge 38$, \ $R(4,6;3)\ge 71$ \ and \ $R(5,5;3)\ge 99$.
\end{theorem}

\begin{theorem}\label{thm:four}
$R(5,6;4)\ge \Bfs$ \ and \ $R(6,6;4)\ge 601$.
\end{theorem}

\begin{theorem}\label{thm:five}
$R(6,6;5)\ge 84$.
\end{theorem}

The proofs consist of explicit colourings. Those for
Theorem~\ref{thm:main}, of orders $37$, $70$ and $98$, are invariant under
groups of affine maps of $\Z_n$ and are described in
Section~\ref{sec:colourings}. Those for Theorem~\ref{thm:four}, of orders
$200$ and $600$, are colourings of the $4$-subsets of a projective line that
are invariant under its group $\mathrm{PSL}(2,q)$ or $\mathrm{PGL}(2,q)$
(Section~\ref{sec:line}). The colouring for
Theorem~\ref{thm:five}, of order $83$, is invariant under all affine maps of
$\Z_{83}$ (Section~\ref{sec:five}). Section~\ref{sec:verification} explains
how the colourings were verified and how a reader can verify them again.
Table~\ref{tab:bounds} compares the bounds.

Two known constructions turn these bounds into further ones
(Section~\ref{sec:product}). The survey notes that $R(5,5;3)\ge 88$ and a
product construction of Bruce, Budden and Hiller~\cite{brbh} give
$R(5,5,5;3)\ge 7570$; the same construction applies to our colouring of
order $98$. And the recurrences of Song, Ye and Liu~\cite{sonyl}, from which
the survey derives its bounds for $R(p,q;4)$, apply to our bound for
$R(5,6;4)$.

\begin{corollary}\label{cor:three}
$R_k(5;3)\ge 98^{2^{k-2}}+1$ for every $k\ge 2$. In particular
$R(5,5,5;3)\ge 9605$.
\end{corollary}

\begin{corollary}\label{cor:four}
$R(5,7;4)\ge \Bfseven$, \ $R(6,7;4)\ge \Bsseven$, \ $R(7,7;4)\ge \Bsevenseven$
\ and \ $R(8,8;4)\ge \Beighteight$.
\end{corollary}

\begin{table}[ht]
\centering
\small
\begin{tabular}{llcc}
\toprule
Number & Previous lower bounds & New lower bound & Order of the colouring\\
\midrule
$R(4,5;3)$ & $35$ \cite{dyb2} & $\mathbf{38}$ & $37$\\
$R(4,6;3)$ & $63$ \cite{dyb3}; \ $64$ \cite{yang} & $\mathbf{71}$ & $70$\\
$R(5,5;3)$ & $88$ \cite{dyb3} & $\mathbf{99}$ & $98$\\
$R(5,6;4)$ & $67$ \cite{survey}; \ $159$ \cite{ozdemir} & $\mathbf{201}$ & $200$\\
$R(6,6;4)$ & $133$ \cite{survey}; \ $317$ \cite{ozdemir} & $\mathbf{601}$ & $600$\\
$R(6,6;5)$ & $72$ \cite{survey,liu} & $\mathbf{84}$ & $83$\\
\midrule
$R(5,5,5;3)$ & $7570$ \cite{survey,brbh} & $\mathbf{9605}$ & $98^2$\\
$R(6,7;4)$ & $661$ \cite{survey}; \ $1581^{\dagger}$ & $\mathbf{2001}$ & --\\
$R(7,7;4)$ & $3961$ \cite{survey}; \ $9481^{\dagger}$ & $\mathbf{12\,001}$ & --\\
$R(8,8;4)$ & $194\,041$ \cite{survey}; \ $464\,521^{\dagger}$ & $\mathbf{588\,001}$ & --\\
\bottomrule
\end{tabular}
\caption{Previous and new lower bounds. The first six rows are
Theorems~\ref{thm:main}, \ref{thm:four} and~\ref{thm:five}; the last four
are Corollaries~\ref{cor:three} and~\ref{cor:four}, obtained from them by
known constructions. The values marked $\dagger$ are not stated
in~\cite{ozdemir}: we derive them from its bound $R(5,6;4)\ge 159$ by the
recurrences of Section~\ref{sec:product}.}
\label{tab:bounds}
\end{table}

\section{The colourings of triples}\label{sec:colourings}

Let $n\ge 2$ and let $a$ be a unit of $\Z_n$ of multiplicative order $e$. We
write
\[
  \G(n,a)=\{\,x\mapsto a^i x+b \;:\; 0\le i<e,\ b\in\Z_n\,\}
\]
for the group of $en$ affine permutations of $\Z_n$ generated by the
translations and the multiplication by $a$. We also let $\G(n,a)$ act on
$\Z_n\cup\{\infty\}$, with every map fixing the point $\infty$. A colouring
of the triples of the point set is \emph{$\G(n,a)$-invariant} if every map of
the group preserves the colour of every triple. Such a colouring is constant
on the orbits of the group on triples, so it is described by naming the orbits
that receive colour~$0$.

We define three colourings.

\begin{itemize}
\item $C_{37}$ is a $\G(37,10)$-invariant colouring of the triples of
  $\Z_{37}$.
\item $C_{70}$ is a $\G(69,4)$-invariant colouring of the triples of
  $\Z_{69}\cup\{\infty\}$.
\item $C_{98}$ is a $\G(97,4)$-invariant colouring of the triples of
  $\Z_{97}\cup\{\infty\}$.
\end{itemize}

In each case colour $0$ is the union of the orbits of the triples listed in
Appendix~\ref{app:reps}, and every other triple has colour~$1$.
Table~\ref{tab:colourings} gives the parameters.

\begin{table}[ht]
\centering
\small
\begin{tabular}{lccc}
\toprule
 & $C_{37}$ & $C_{70}$ & $C_{98}$\\
\midrule
Points & $\Z_{37}$ & $\Z_{69}\cup\{\infty\}$ & $\Z_{97}\cup\{\infty\}$\\
Multiplier $a$ (its order) & $10$ ($3$) & $4$ ($11$) & $4$ ($24$)\\
Order of the group $\G(n,a)$ & $111$ & $759$ & $2328$\\
\addlinespace
Orbits of triples & $78$ & $74$ & $72$\\
\quad by size & $66\cdot 111+12\cdot 37$ & $72\cdot 759+69+23$ & $60\cdot 2328+8\cdot 1164+4\cdot 776$\\
Orbits in colour $0$ & $32$ & $27$ & $36$\\
\quad by size & $27\cdot 111+5\cdot 37$ & $26\cdot 759+69$ & $30\cdot 2328+4\cdot 1164+2\cdot 776$\\
\addlinespace
Triples in colour $0$ & $3182$ & $19\,803$ & $76\,048$\\
Triples in colour $1$ & $4588$ & $34\,937$ & $76\,048$\\
\bottomrule
\end{tabular}
\caption{The three colourings of triples. A term $m\cdot d$ stands for $m$
orbits of size $d$.}
\label{tab:colourings}
\end{table}

\begin{proposition}\label{prop:valid}
$C_{37}$ is a $(4,5;3)$-colouring of order $37$, \ $C_{70}$ is a
$(4,6;3)$-colouring of order $70$, \ and $C_{98}$ is a $(5,5;3)$-colouring of
order $98$.
\end{proposition}

\begin{proof}[Computational verification]
By the exhaustive search described in Section~\ref{sec:verification}.
\end{proof}

Theorem~\ref{thm:main} follows at once from
Proposition~\ref{prop:valid}.

\begin{remark}[a colour-reversing symmetry]\label{rem:reverse}
The map $x\mapsto 2x$ of $\Z_{97}\cup\{\infty\}$, which fixes $\infty$,
exchanges the two colours of $C_{98}$ on every triple. Since $2$ generates the
squares modulo $97$ and $4$ generates the fourth powers, the multiplications
by fourth powers preserve $C_{98}$ and the multiplications by the other
squares reverse it. In particular the two colour classes of $C_{98}$ are
isomorphic, and each has $\binom{98}{3}/2=76\,048$ triples.
\end{remark}

\begin{remark}\label{rem:notprojective}
$C_{98}$ is not invariant under the projective map $x\mapsto -1/x$ of the
line $\Z_{97}\cup\{\infty\}$: that map preserves the colour of only $74\,512$
of the $152\,096$ triples. So $C_{98}$ is not invariant under
$\mathrm{PSL}(2,97)$ in its action on this line.
\end{remark}

\section{Colourings of the projective line}\label{sec:line}

Let $q$ be a prime and let $L_q=\mathbb{F}_q\cup\{\infty\}$ be the projective
line over the field with $q$ elements, a set of $q+1$ points. The group
$\mathrm{PGL}(2,q)$ acts on $L_q$ by the maps
$x\mapsto(\alpha x+\beta)/(\gamma x+\delta)$ with
$\alpha\delta-\beta\gamma\ne 0$. For four distinct points $a,b,c,d$ of $L_q$
the \emph{cross-ratio} is
\[
  [a,b;c,d]=\frac{(a-c)(b-d)}{(a-d)(b-c)},
\]
where the two factors that contain $\infty$, if any, are cancelled against
each other; for instance $[a,b;c,\infty]=(a-c)/(b-c)$. The cross-ratio lies in
$\mathbb{F}_q\setminus\{0,1\}$ and is unchanged by every map of
$\mathrm{PGL}(2,q)$. Reordering the four points replaces a value $\lambda$ by
one of
\[
  \lambda,\qquad \frac1\lambda,\qquad 1-\lambda,\qquad \frac1{1-\lambda},\qquad
  \frac{\lambda}{\lambda-1},\qquad \frac{\lambda-1}{\lambda},
\]
and we call this set of at most six values the \emph{class} of $\lambda$.
Since $\mathrm{PGL}(2,q)$ is sharply $3$-transitive on $L_q$, two $4$-subsets
of $L_q$ lie in one orbit of the group exactly when their cross-ratios lie in
one class. Hence a $\mathrm{PGL}(2,q)$-invariant colouring of the $4$-subsets
of $L_q$ is the same thing as a function
$\tau\colon\mathbb{F}_q\setminus\{0,1\}\to\{0,1\}$ that is constant on
classes: the colour of $\{a,b,c,d\}$ is $\tau([a,b;c,d])$. There are about
$q/6$ classes, so such a colouring is a short table.

Colourings of the $4$-subsets of $L_q$ that are invariant under the subgroup
$\mathrm{PSL}(2,q)$ were used by \"Ozdemir~\cite{ozdemir}, with $q=47$ for
$R(5,5;4)$ and $q=157$ for $R(5,6;4)$. We searched the
$\mathrm{PGL}(2,q)$-invariant colourings, for $(5,6;4)$ at larger $q$ and for
$(6,6;4)$ directly, and define three colourings by their functions $\tau$.

\begin{itemize}
\item $D_{194}$ is a colouring of the $4$-subsets of $L_{193}$; its $\tau$
  takes the value $0$ on $15$ of the $33$ classes.
\item $E_{500}$ is a colouring of the $4$-subsets of $L_{499}$; its $\tau$
  takes the value $0$ on $41$ of the $84$ classes.
\item $E_{600}$ is a colouring of the $4$-subsets of $L_{599}$; its $\tau$
  takes the value $0$ on $51$ of the $100$ classes.
\end{itemize}

Appendix~\ref{app:tables} lists these classes. $D_{194}$ has $26\,059\,632$
of its $57\,211\,376$ $4$-subsets in colour $0$; $E_{500}$ has $1\,273\,572\,750$ of
$2\,573\,031\,125$; $E_{600}$ has $2\,740\,245\,300$ of $5\,346\,164\,850$.

\begin{proposition}\label{prop:line}
$D_{194}$ is a $(5,6;4)$-colouring of order $194$. \ $E_{500}$ and $E_{600}$
are $(6,6;4)$-colourings of orders $500$ and $600$.
\end{proposition}

\begin{proof}[Computational verification]
By the exhaustive searches described in Section~\ref{sec:verification}.
\end{proof}

The second bound of Theorem~\ref{thm:four} follows from $E_{600}$; the first is
proved below, with a smaller group. We record $E_{500}$ as well because its
verification uses no symmetry at all.

Appendix~\ref{app:census} records the primes $q$ at which our searches found
colourings of these kinds, and those at which the solver reported that the
family contains none.

\paragraph{The smaller group.}
The subgroup $\mathrm{PSL}(2,q)$ of $\mathrm{PGL}(2,q)$ consists of the maps
$x\mapsto(\alpha x+\beta)/(\gamma x+\delta)$ for which
$\alpha\delta-\beta\gamma$ is a nonzero square. It is the group used
in~\cite{ozdemir}. Each orbit of $\mathrm{PGL}(2,q)$ on the $4$-subsets of
$L_q$ is one orbit of the subgroup or splits into two, so the subgroup admits
more invariant colourings. Let
$q\equiv 3\pmod 4$, so that $-1$ is not a square. For distinct
$a,b,c\in\mathbb{F}_q$ the map
\[
  g\colon x\mapsto\frac{(c-a)(x-b)}{(c-b)(x-a)}
\]
sends $a$, $b$, $c$ to $\infty$, $0$, $1$ and any fourth point $d$ to
$[a,b;c,d]$, and its determinant is
\[
  \Delta(a,b,c)=(b-a)(c-a)(c-b).
\]
If $\Delta(a,b,c)$ is a square, $g$ lies in $\mathrm{PSL}(2,q)$. If not, the
map $x\mapsto 1/g(x)$ does, and it sends $a$, $b$, $c$, $d$ to $0$, $\infty$,
$1$ and $1/[a,b;c,d]$. So every $4$-subset of $L_q$ lies in the orbit of a
set $\{\infty,0,1,\lambda\}$, and a $\mathrm{PSL}(2,q)$-invariant colouring is
known from the function $\sigma$ on $\mathbb{F}_q\setminus\{0,1\}$ for which
$\sigma(\lambda)$ is the colour of $\{\infty,0,1,\lambda\}$: for distinct
$a,b,c\in\mathbb{F}_q$ and a fourth point $d$ of $L_q$,
\begin{equation}\label{eq:psl}
  \text{the colour of }\{a,b,c,d\}\text{ is }
  \begin{cases}
    \sigma\bigl([a,b;c,d]\bigr) & \text{if $\Delta(a,b,c)$ is a square,}\\[2pt]
    \sigma\bigl(1/[a,b;c,d]\bigr) & \text{otherwise.}
  \end{cases}
\end{equation}
The next lemma says which functions $\sigma$ occur. We call the set of
values that it describes the \emph{narrow class} of $\lambda$.

\begin{lemma}\label{lem:narrow}
Let $q\equiv 3\pmod 4$ and $\lambda\in\mathbb{F}_q\setminus\{0,1\}$. The
values $\mu$ for which $\{\infty,0,1,\mu\}$ lies in the orbit of
$\{\infty,0,1,\lambda\}$ under $\mathrm{PSL}(2,q)$ are $\lambda$,
$1/(1-\lambda)$ and $(\lambda-1)/\lambda$ if $-\lambda$ and $\lambda-1$ are
both squares, and all the values in the class of $\lambda$ otherwise.
\end{lemma}

\begin{proof}
Write $S=\{\infty,0,1,\lambda\}$. The map $x\mapsto 1/(1-x)$ has determinant
$1$ and permutes $\infty$, $0$, $1$ cyclically; it and its square send
$\lambda$ to $1/(1-\lambda)$ and to $(\lambda-1)/\lambda$, so these values
always occur. The maps $x\mapsto\lambda/x$ and $x\mapsto(x-\lambda)/(x-1)$
send $S$ to itself and have determinants $-\lambda$ and $\lambda-1$; with
their product and the identity they form a group $V$ that permutes the four
points of $S$ transitively.

Suppose that $-\lambda$ or $\lambda-1$ is not a square, so that some
$v\in V$ lies outside $\mathrm{PSL}(2,q)$. If $\mu$ lies in the class of
$\lambda$, some $h\in\mathrm{PGL}(2,q)$ sends $S$ to $\{\infty,0,1,\mu\}$,
and so does $hv$; one of the two lies in $\mathrm{PSL}(2,q)$. Suppose now
that both are squares, so that $V$ lies in $\mathrm{PSL}(2,q)$, and let
$h\in\mathrm{PSL}(2,q)$ send $S$ to $\{\infty,0,1,\mu\}$. Choosing $v\in V$
with $hv(\lambda)=\mu$, we obtain a map $hv$ of $\mathrm{PSL}(2,q)$ that
sends $\{\infty,0,1\}$ to itself. The six maps of $\mathrm{PGL}(2,q)$ that do
this are $x$, $1/(1-x)$ and $(x-1)/x$, of determinant $1$, and $1/x$, $1-x$
and $x/(x-1)$, of determinant $-1$, which is not a square. So $hv$ is one of
the first three, and $\mu=hv(\lambda)$ is one of $\lambda$, $1/(1-\lambda)$
and $(\lambda-1)/\lambda$.
\end{proof}

Hence, for $q\equiv 3\pmod 4$, the $\mathrm{PSL}(2,q)$-invariant colourings
of the $4$-subsets of $L_q$ are given through~\eqref{eq:psl} by the functions
$\sigma$ that are constant on narrow classes, and the orbits of the group on
$4$-subsets correspond to the narrow classes.

We define $F_{200}$ as the $\mathrm{PSL}(2,199)$-invariant colouring of the
$4$-subsets of $L_{199}$ whose function $\sigma$ takes the value $0$ on the
$20$ narrow classes listed in Appendix~\ref{app:tables}, out of $43$, and
the value $1$ on the others. It has $27\,909\,750$ of its $64\,684\,950$
$4$-subsets in colour $0$.

\begin{proposition}\label{prop:psl}
$F_{200}$ is a $(5,6;4)$-colouring of order $200$.
\end{proposition}

\begin{proof}[Computational verification]
By the exhaustive search described in Section~\ref{sec:verification}.
\end{proof}

So $R(5,6;4)\ge 201$, which is the first bound of Theorem~\ref{thm:four}.
The colouring $D_{194}$ gives $R(5,6;4)\ge 195$ with the larger group.

\section{A colouring of \texorpdfstring{$5$}{5}-subsets}\label{sec:five}

Colourings with the full symmetry of a projective line cannot be used for
$R(6,6;5)$, for the following reason.

\begin{lemma}\label{lem:obstruction}
Let a group $G$ act on a set $X$, let $t\ge 2$, and let $S$ be a set of $t+1$
points of $X$. If some subgroup of $G$ maps $S$ to itself and permutes its
points transitively, then in every $G$-invariant colouring of the $t$-subsets
of $X$ all $t$-subsets of $S$ have the same colour. In particular there is
then no $G$-invariant $(t+1,t+1;t)$-colouring of $X$.
\end{lemma}

\begin{proof}
The $t$-subsets of $S$ are the complements in $S$ of its points, so the
subgroup permutes them transitively, and an invariant colouring gives them one
colour.
\end{proof}

For every prime $q\ge 5$ the lemma applies to $\mathrm{PGL}(2,q)$ acting on
$L_q$, with $t=5$. For $q=5$ the group permutes the six points of $L_5$
transitively. For $q=7$ the maps $x\mapsto ax$ with $a\ne 0$ permute the six
points $1,\dots,6$ transitively. For $q\ge 11$ consider the subgroup that maps
$\{0,1,\infty\}$ to itself. It has six elements. Each of its three elements
of order two fixes at most two points of $L_q$, and its two elements of order
three, which are inverse to each other, fix the same points, at most two. So
at most eight points are fixed by a non-identity element, and some point of
$L_q$ has an orbit of six points, on which the subgroup acts transitively. In
each case Lemma~\ref{lem:obstruction} shows that no
$\mathrm{PGL}(2,q)$-invariant $(6,6;5)$-colouring of $L_q$ exists.

Let $A_p$ be the group of all maps $x\mapsto ax+b$ of $\Z_p$ with $a\ne 0$,
for a prime $p\ge 5$. It has $p(p-1)$ elements. If $p\equiv 1\pmod 6$, the
lemma applies to $A_p$ as well: the six solutions of $a^6=1$ in $\Z_p$ are
permuted transitively by the multiplications by these same six numbers. If
$p\equiv 5\pmod 6$, it does not apply: a group that permutes six points
transitively has an order divisible by~$6$, and $3$ does not divide
$p(p-1)$. This is the case for $p=83$, where $p(p-1)=2\cdot 41\cdot 83$.

The group $A_p$ is sharply $2$-transitive: for two points $u\ne v$ the map
$x\mapsto(x-u)/(v-u)$ sends $u$ to $0$ and $v$ to $1$. So every $5$-subset of
$\Z_p$ can be moved to one that contains $0$ and $1$, and an $A_p$-invariant
colouring of the $5$-subsets is determined by its colours on the
$\binom{p-2}{3}$ subsets that contain $0$ and~$1$.

We define $H_{83}$ as the $A_{83}$-invariant colouring of the $5$-subsets of
$\Z_{83}$ whose colours on the $\binom{81}{3}=85\,320$ subsets through $0$ and
$1$ are given in the ancillary file \path{anc/R665_83_points_short.json}. It
has $14\,500\,183$ of its $29\,034\,396$ $5$-subsets in colour $0$, and the $5$-subsets
fall into $4276$ orbits of $A_{83}$.

\begin{proposition}\label{prop:five}
$H_{83}$ is a $(6,6;5)$-colouring of order $83$.
\end{proposition}

\begin{proof}[Computational verification]
By the exhaustive search described in Section~\ref{sec:verification}.
\end{proof}

Theorem~\ref{thm:five} follows.

\section{Verification}\label{sec:verification}

\paragraph{The check.}
Let $t$ be the size of the subsets that are coloured, $c$ a colour and $s$ a
size. A set of points all of whose $t$-subsets have colour $c$ is called
\emph{one-coloured}. Our programs search for a one-coloured set of $s$
points by backtracking: a one-coloured set is extended by a point larger than
all of its points, and the new point is accepted only if every $t$-subset
that consists of it and $t-1$ points of the set has colour $c$. The search is
exhaustive, because every initial segment of a one-coloured set, in
increasing order of its points, is one-coloured. A colouring is accepted as
a $(k,l;t)$-colouring only if the search finds no one-coloured set of $k$
points in colour $0$ and none of $l$ points in colour $1$. For the
colourings of triples the sets so excluded number $\binom{37}{4}=66\,045$ and
$\binom{37}{5}=435\,897$ for $C_{37}$; $\binom{70}{4}=916\,895$ and
$\binom{70}{6}=131\,115\,985$ for $C_{70}$; and
$\binom{98}{5}=67\,910\,864$ for each colour of $C_{98}$.

\paragraph{Independent programs for the colourings of triples.}
Each of $C_{37}$, $C_{70}$ and $C_{98}$ was checked by several independently written programs, none of
which shares code with the search that produced the colourings. One of them,
\path{anc/verify.py}, accompanies this paper. It uses only the Python
standard library and reads only a plain list of the triples with their
colours, so it depends on no ordering convention. It was also run on altered colourings in which a one-coloured set had been
planted, and reported it. On an ordinary laptop \path{anc/verify.py} needs a few seconds for each
colouring.

\paragraph{From the paper to the files.}
The script \path{anc/rebuild_from_rule.py} rebuilds each colouring from its
group and the representatives of Appendix~\ref{app:reps}, and compares the
result with the stored colouring. So the printed data and the ancillary files
describe the same three objects. Section~\ref{sec:data} gives SHA-256 digests.

\paragraph{The colourings of $4$-subsets.}
For $D_{194}$, $E_{500}$ and $E_{600}$ the colour of a $4$-subset is computed
from the table $\tau$ and the cross-ratio of its points taken in increasing
order, with $\infty$ as the last point. The program
\path{anc/verify_crossratio.c} treats this as a plain function on
$4$-subsets and runs the search described above, with bit vectors, for
$4$-subsets in place of triples. It assumes nothing about the table and uses
no symmetry.
\begin{itemize}
\item For $D_{194}$ it examined all $83\,870\,080$ one-coloured sets of $5$
  points in colour $1$; none extends to $6$ points, and no $5$ points are
  one-coloured in colour $0$.
\item For $E_{500}$ it examined $6\,411\,351\,600$ one-coloured sets of $5$ points
  in colour $0$ and $6\,647\,428\,500$ in colour $1$; none extends to $6$ points.
\item For $E_{600}$ the program first confirmed, on each of the $5\,346\,164\,850$
  $4$-subsets, that the map $x\mapsto x+1$ of $L_{599}$, which fixes
  $\infty$, preserves the colour. Every set of at least two points contains a
  point other than $\infty$ and so has a translate through the point $0$;
  the program then examined the $147\,041\,921$ and $121\,788\,680$ one-coloured
  sets of $5$ points through $0$ in the two colours, and none extends to $6$
  points.
\end{itemize}
Each run began with controls: a colouring known to be valid had to pass, and
tables altered on purpose had to fail, among them one that $x\mapsto x+1$
does not preserve, which the program must refuse. On $E_{318}$, the analogous
colouring for $q=317$, the search through $0$ and the search without symmetry
gave the same answer.

$D_{194}$ was checked a second time in a different way: the colouring was
written out in full, one colour for each of its $57\,211\,376$ $4$-subsets, and
the program \path{anc/verify_rank.c} examined every one of the
$\binom{194}{5}=2\,174\,032\,288$ sets of $5$ points and of the
$\binom{194}{6}=68\,482\,017\,072$ sets of $6$ points. The two programs also agree
on the number of $4$-subsets in colour $0$.

\paragraph{A second program for the tables $\tau$.}
The program \path{anc/verify_projective.cpp} was written independently of
the first. It checks that a table is constant on the classes of cross-ratio
values, so that the colouring is invariant under $\mathrm{PGL}(2,q)$, and then
examines every set of the size in question that contains $0$, $1$ and
$\infty$. This is a complete check, because $\mathrm{PGL}(2,q)$ is
$3$-transitive on $L_q$: every one-coloured set of at least three points has
a one-coloured image that contains these three points. For $D_{194}$ the
program examined $\binom{191}{2}=18\,145$ sets of $5$ points and
$\binom{191}{3}=1\,143\,135$ sets of $6$ points; for $E_{500}$,
$\binom{497}{3}=20\,337\,240$ sets of $6$ points; and for $E_{600}$,
$\binom{597}{3}=35\,284\,690$ sets of $6$ points. None of them is
one-coloured in a forbidden colour. The script
\path{anc/check_projective.py} compiles the program and runs it on the three
tables.

\paragraph{The colouring $F_{200}$.}
$F_{200}$ was written out in full, one colour for each of its $64\,684\,950$
$4$-subsets, and \path{anc/verify_rank.c} examined every one of the
$\binom{200}{5}=2\,535\,650\,040$ sets of $5$ points and of the
$\binom{200}{6}=82\,408\,626\,300$ sets of $6$ points. This check uses no
symmetry and no property of the group.
The script \path{anc/rebuild_psl_rule.py} confirms that the function
$\sigma$ is constant on narrow classes, rebuilds all the colours from
$\sigma$ and~\eqref{eq:psl}, and compares the SHA-256 digest of the result
with that of the colouring written out in full; so the description in
Section~\ref{sec:line} and the file that was checked are the same colouring.
It also examines the $\binom{197}{2}=19\,306$ sets of $5$ points and the
$\binom{197}{3}=1\,254\,890$ sets of $6$ points that contain $\infty$, $0$
and $1$, and finds none that is one-coloured in a forbidden colour. Every set
of at least four points of $L_{199}$ has an image under
$\mathrm{PSL}(2,199)$ that contains these three points, so this is a second
complete check, by another program and another method.

\paragraph{The colouring of $5$-subsets.}
$H_{83}$ was written out in full from its short form, and
\path{anc/verify_rank.c} examined every one of the $\binom{83}{6}=377\,447\,148$
sets of $6$ points, for each colour. The search program's own check, written
independently of it and using no symmetry, agrees. The script
\path{anc/rebuild_short_form.py} rebuilds the full colouring from the short
form and compares its SHA-256 digest.

\section{Consequences}\label{sec:product}

The following construction is the one of Bruce, Budden and
Hiller~\cite{brbh} that the survey~\cite{survey} uses; we include the short
argument for completeness.

\begin{lemma}\label{lem:product}
Let $k\ge 2$. Suppose that the triples of an $N$-set $X$ can be coloured with
$k$ colours so that no $5$ points have all their triples in one colour. Then
the triples of $X\times X$ can be coloured with $k+1$ colours with the same
property.
\end{lemma}

\begin{proof}
Let $\chi$ be the given colouring, with colours $1,\dots,k$. Colour a triple
\[
  T=\{(x_1,y_1),(x_2,y_2),(x_3,y_3)\}
\]
of $X\times X$ as follows. If
$x_1,x_2,x_3$ are distinct, $T$ gets the colour $\chi(\{x_1,x_2,x_3\})$. If
$x_1=x_2=x_3$, then $y_1,y_2,y_3$ are distinct and $T$ gets the colour
$\chi(\{y_1,y_2,y_3\})$. Otherwise exactly two of $x_1,x_2,x_3$ are equal, and
$T$ gets the new colour $k+1$.

Let $S$ be a set of $5$ points of $X\times X$ whose triples all have the same
colour $c$. Suppose first that $c\le k$. If two points of $S$ have the same
first coordinate and a third point has a different one, these three points
form a triple of colour $k+1$. So the first coordinates of the points of $S$
are either all distinct or all equal. In the first case they form $5$ points
of $X$ whose triples all have colour $c$ under $\chi$; in the second case the
second coordinates do. Both contradict the assumption on $\chi$. Suppose now
that $c=k+1$. If three points of $S$ have distinct first coordinates, their
triple has a colour at most $k$. Otherwise the $5$ points of $S$ have at most
two distinct first coordinates, so three of them share one, and their triple
has a colour at most $k$. In both cases we have a contradiction.
\end{proof}

\begin{proof}[Proof of Corollary~\ref{cor:three}]
By Proposition~\ref{prop:valid}, the triples of a $98$-set can be
coloured with $2$ colours so that no $5$ points have all their triples in one
colour. Applying Lemma~\ref{lem:product} $k-2$ times gives such a colouring
with $k$ colours on $98^{2^{k-2}}$ points. For $k=3$ this is $98^2=9604$
points.
\end{proof}

\begin{proof}[Proof of Corollary~\ref{cor:four}]
The survey~\cite[Section~7.1]{survey} quotes two recurrences of Song, Ye and
Liu~\cite{sonyl}. We use them only in the following ranges:
\begin{align*}
  R(p,q;4)&\ge 2R(p-1,q;4)-1 &&\text{for } p\ge 6,\ q\ge 5,\\
  R(p,q;4)&\ge (p-1)R(p-1,q;4)-p+2 &&\text{for } p\ge 6,\ q\ge 7.
\end{align*}
Together with the symmetry $R(p,q;4)=R(q,p;4)$ and Theorem~\ref{thm:four},
the first gives $R(7,5;4)\ge 2\cdot 201-1=401$. The second, applied to
$(p,q)=(6,7)$, $(7,7)$, $(8,7)$ and $(8,8)$ in turn, gives
$R(6,7;4)\ge 5\cdot 401-4=2001$, \ $R(7,7;4)\ge 6\cdot 2001-5=12\,001$, \
$R(8,7;4)\ge 7\cdot 12\,001-6=84\,001$ \ and \
$R(8,8;4)\ge 7\cdot 84\,001-6=588\,001$. The values $661$, $3961$ and
$194\,041$ that the survey lists for $R(6,7;4)$, $R(7,7;4)$ and $R(8,8;4)$
arise in the same way from $R(5,6;4)\ge 67$. From the bound
$R(5,6;4)\ge 159$ of~\cite{ozdemir} the same steps give $1581$, $9481$ and
$464\,521$; these are the values marked $\dagger$ in Table~\ref{tab:bounds}.
\end{proof}

\section{Remarks on the search}\label{sec:search}

\paragraph{Method.}
The colourings of triples were found by exact search over invariant colourings. For a
chosen group $\G(n,a)$, with or without a fixed point, there is one Boolean
unknown for each orbit of triples. The group contains all translations, and
every set of at least three points contains a point of $\Z_n$, so it can be
translated to a set that contains the point $0$; hence it suffices to forbid
the one-coloured sets of the relevant sizes that contain $0$. The resulting formulas were given to the SAT solvers
Glucose~\cite{glucose} and CaDiCaL~\cite{cadical} through the PySAT
toolkit~\cite{pysat}. For $C_{98}$ we required in addition that $x\mapsto 2x$
exchange the colours (Remark~\ref{rem:reverse}), which leaves $36$ unknowns;
the solver answered in $14$ seconds. $C_{70}$ has $74$ unknowns and took six
minutes.

The tables $\tau$ of Section~\ref{sec:line} were found in the same way, with
one unknown for each class of cross-ratio values: since $\mathrm{PGL}(2,q)$
is $3$-transitive, it suffices to forbid the one-coloured sets that contain
$0$, $1$ and $\infty$. For $F_{200}$ there is one unknown for each narrow
class, and the sets forbidden are those that contain $\infty$, $0$ and one of
two further fixed points, a square and a non-square. For $q\le 499$ the SAT
solver answered at once. For $E_{600}$ and for $H_{83}$, which has $4276$
unknowns, we used a local search of the probSAT type~\cite{probsat} instead;
it found $E_{600}$ in about $6\cdot 10^4$ steps and $H_{83}$ in about
$3\cdot 10^8$. A local search can only find a colouring; what it finds is
checked as in Section~\ref{sec:verification}.

\paragraph{What the search did not find.}
The following are reports of a SAT solver on the encoded families, or of a
local search where that is said; no proof certificate was recorded. They describe the limits of these particular
families, establish no upper bound, and say nothing about other colourings.
\begin{itemize}
\item For $R(5,5;3)$: on $\Z_{97}$ with two fixed points ($99$ points), the
  solver returned UNSAT for the $\G(97,4)$-invariant family and for the
  $\G(97,2)$-invariant family. For the larger $\G(97,16)$-invariant family
  the search ended without an answer.
\item For $R(4,6;3)$: on $\Z_{69}$ with two fixed points ($71$ points), the
  solver returned UNSAT for the $\G(69,4)$-invariant family.
\item For $R(4,5;3)$: on $\Z_{37}$ with one fixed point ($38$ points), the
  solver returned UNSAT for the $\G(37,10)$-invariant family. The searches
  we ran on $\Z_{39}$ and $\Z_{43}$ ended without an answer.
\item For $R(6,6;4)$: on $L_{701}$ ($702$ points), neither the local search
  nor the SAT solver gave an answer in the time allowed.
\item For $R(5,6;4)$: the solver returned UNSAT for the
  $\mathrm{PGL}(2,q)$-invariant family at every prime $q$ with
  $197\le q\le 331$, and for the $\mathrm{PSL}(2,q)$-invariant family at
  every prime $q$ with $211\le q\le 331$ (Appendix~\ref{app:census}).
\item For $R(6,6;5)$: on $\Z_{83}$ with one fixed point ($84$ points), both with all
  multipliers and with the squares as multipliers, and on $\Z_{89}$ with all
  multipliers ($89$ points), the local search found no colouring. A local
  search cannot show that none exists.
\end{itemize}

\section{Data availability}\label{sec:data}

The colourings accompany this paper as ancillary files, which are listed in
Table~\ref{tab:files}.

\begin{itemize}
\item \emph{Triples.} Each of $C_{37}$, $C_{70}$ and $C_{98}$ comes in two
  forms: a JSON file that holds the colouring as a string of the characters
  \texttt{0} and \texttt{1}, one for each triple of the points $0,\dots,n-1$
  in lexicographic order, together with its group and orbit representatives;
  and a CSV file with one row for each triple, giving its three points and
  its colour. The point $\infty$ is the last point, numbered $69$ in $C_{70}$
  and $97$ in $C_{98}$.
\item \emph{$4$-subsets.} Each of $D_{194}$, $E_{500}$ and $E_{600}$ is a
  JSON file with its table $\tau$ as a string of $q$ characters (the places
  $0$ and $1$ are not used) and the list of classes of colour $0$. The points
  are $0,\dots,q-1$ and $\infty$, which is the last point, numbered $q$.
  $F_{200}$ comes in two forms: a JSON file of the same kind, with its table
  $\sigma$ and the list of narrow classes of colour $0$; and a text file,
  compressed with gzip, with one character for each $4$-subset of the points
  $0,\dots,199$ in lexicographic order, the point $\infty$ being numbered
  $199$.
\item \emph{$5$-subsets.} $H_{83}$ is a JSON file with the colours of the
  $85\,320$ subsets through $0$ and $1$, in lexicographic order.
\end{itemize}

\begin{table}[t]
\centering
\small
\begin{tabular}{ll}
\toprule
File & Content\\
\midrule
\path{anc/R453_37_points.json}, \path{.csv} & $C_{37}$\\
\path{anc/R463_70_points.json}, \path{.csv} & $C_{70}$\\
\path{anc/R553_98_points.json}, \path{.csv} & $C_{98}$\\
\path{anc/R564_194_points_rule.json} & $D_{194}$\\
\path{anc/R564_200_points_rule.json} & $F_{200}$: the table $\sigma$\\
\path{anc/R564_200_points.txt.gz} & $F_{200}$, written out in full\\
\path{anc/R664_500_points_rule.json} & $E_{500}$\\
\path{anc/R664_600_points_rule.json} & $E_{600}$\\
\path{anc/R665_83_points_short.json} & $H_{83}$\\
\path{anc/verify.py} & the check of Section~\ref{sec:verification}, for triples\\
\path{anc/rebuild_from_rule.py} & rebuilds a colouring from Appendix~\ref{app:reps}\\
\path{anc/verify_crossratio.c} & the check for $D_{194}$, $E_{500}$, $E_{600}$\\
\path{anc/verify_projective.cpp} & the second check for the same three\\
\path{anc/check_projective.py} & compiles and runs it\\
\path{anc/verify_rank.c} & examines every set, for a colouring written out in full\\
\path{anc/rebuild_short_form.py} & writes $H_{83}$ out in full\\
\path{anc/rebuild_psl_rule.py} & rebuilds $F_{200}$ from $\sigma$ and compares\\
\path{anc/README.md} & formats and usage\\
\bottomrule
\end{tabular}
\caption{The ancillary files.}
\label{tab:files}
\end{table}

\noindent SHA-256 digests of the colouring strings of $C_{37}$, $C_{70}$,
$C_{98}$, $F_{200}$ and $H_{83}$, and of the tables $\tau$ of $D_{194}$, $E_{500}$ and
$E_{600}$:
{\small
\begin{itemize}
\item $C_{37}$: \texttt{995fb6bb93ca2b0be5809d981cefbaa7\allowbreak 2712ce8b1463809dd66376af9f6c0e16}
\item $C_{70}$: \texttt{f7c37ea5e0b72ad3773b13c0a291697f\allowbreak 336998170931f208f6588352ce31c662}
\item $C_{98}$: \texttt{53e042c9f44d870cdccd4d86b5e45003\allowbreak fadf981d7757580f863632a945f5e17d}
\item $D_{194}$: \texttt{477e5506f17a2db541c806f9c57ca20c\allowbreak 2e1d9adde4918dd01c29075d3092dbe7}
\item $F_{200}$: \texttt{621a303b0e947fb3e52799288ba7f8ef\allowbreak 883f3260f6a7e260c843eb902001602d}
\item $E_{500}$: \texttt{3ed928e24dabb537c6d26c5da0251af9\allowbreak 562e40f9be4c9a2f82b7ff55a97931a1}
\item $E_{600}$: \texttt{c4f4b30a65c395b1d7e8fd014a6cd17a\allowbreak 31acb70672ef7294da95a934cf4a7a6f}
\item $H_{83}$: \texttt{b196450e13f22883d9e9c974c3f203b0\allowbreak a958db0560110c2cde73c5bd803daadb}
\end{itemize}}

\subsection*{Acknowledgments}
\phantomsection\addcontentsline{toc}{subsection}{Acknowledgments}
The author thanks Stanis{\l}aw P.\ Radziszowski for maintaining the dynamic
survey on which this work relies. Computational searches were conducted using
the author's Solon research software. AI-based tools were used in the searches
and in the preparation of this manuscript; the author is responsible for all
of its content. The reported lower bounds are independently verifiable from
the ancillary files.

\appendix
\section{Orbit representatives}\label{app:reps}

For each colouring we list one triple from each orbit in colour $0$; the
triple listed is the lexicographically first of its orbit. Colour $0$ is the
union of the orbits of these triples under the group named, and every other
triple has colour $1$.

\paragraph{$C_{37}$: $\G(37,10)$ on $\Z_{37}$, $32$ orbits.}
{\small\raggedright
$\{0,1,2\}$,\allowbreak\ $\{0,1,3\}$,\allowbreak\ $\{0,1,4\}$,\allowbreak\ $\{0,1,5\}$,\allowbreak\ $\{0,1,6\}$,\allowbreak\ $\{0,1,12\}$,\allowbreak\ $\{0,1,14\}$,\allowbreak\ $\{0,1,20\}$,\allowbreak\ $\{0,1,21\}$,\allowbreak\ $\{0,1,24\}$,\allowbreak\ $\{0,1,25\}$,\allowbreak\ $\{0,1,30\}$,\allowbreak\ $\{0,1,33\}$,\allowbreak\ $\{0,2,9\}$,\allowbreak\ $\{0,2,14\}$,\allowbreak\ $\{0,2,15\}$,\allowbreak\ $\{0,2,18\}$,\allowbreak\ $\{0,2,22\}$,\allowbreak\ $\{0,2,23\}$,\allowbreak\ $\{0,2,31\}$,\allowbreak\ $\{0,2,32\}$,\allowbreak\ $\{0,2,33\}$,\allowbreak\ $\{0,3,6\}$,\allowbreak\ $\{0,3,8\}$,\allowbreak\ $\{0,3,9\}$,\allowbreak\ $\{0,3,16\}$,\allowbreak\ $\{0,3,19\}$,\allowbreak\ $\{0,5,18\}$,\allowbreak\ $\{0,5,21\}$,\allowbreak\ $\{0,6,14\}$,\allowbreak\ $\{0,6,29\}$,\allowbreak\ $\{0,9,25\}$.
\par}

\paragraph{$C_{70}$: $\G(69,4)$ on $\Z_{69}\cup\{\infty\}$, $27$ orbits.}
{\small\raggedright
$\{0,1,3\}$,\allowbreak\ $\{0,1,5\}$,\allowbreak\ $\{0,1,8\}$,\allowbreak\ $\{0,1,9\}$,\allowbreak\ $\{0,1,10\}$,\allowbreak\ $\{0,1,13\}$,\allowbreak\ $\{0,1,18\}$,\allowbreak\ $\{0,1,20\}$,\allowbreak\ $\{0,1,23\}$,\allowbreak\ $\{0,1,28\}$,\allowbreak\ $\{0,1,33\}$,\allowbreak\ $\{0,1,36\}$,\allowbreak\ $\{0,1,37\}$,\allowbreak\ $\{0,1,42\}$,\allowbreak\ $\{0,1,43\}$,\allowbreak\ $\{0,1,44\}$,\allowbreak\ $\{0,1,55\}$,\allowbreak\ $\{0,1,60\}$,\allowbreak\ $\{0,1,62\}$,\allowbreak\ $\{0,2,39\}$,\allowbreak\ $\{0,2,47\}$,\allowbreak\ $\{0,2,48\}$,\allowbreak\ $\{0,2,\infty\}$,\allowbreak\ $\{0,3,15\}$,\allowbreak\ $\{0,3,21\}$,\allowbreak\ $\{0,3,27\}$,\allowbreak\ $\{0,23,\infty\}$.
\par}

\paragraph{$C_{98}$: $\G(97,4)$ on $\Z_{97}\cup\{\infty\}$, $36$ orbits.}
{\small\raggedright
$\{0,1,7\}$,\allowbreak\ $\{0,1,10\}$,\allowbreak\ $\{0,1,11\}$,\allowbreak\ $\{0,1,12\}$,\allowbreak\ $\{0,1,13\}$,\allowbreak\ $\{0,1,15\}$,\allowbreak\ $\{0,1,19\}$,\allowbreak\ $\{0,1,20\}$,\allowbreak\ $\{0,1,27\}$,\allowbreak\ $\{0,1,29\}$,\allowbreak\ $\{0,1,30\}$,\allowbreak\ $\{0,1,31\}$,\allowbreak\ $\{0,1,32\}$,\allowbreak\ $\{0,1,33\}$,\allowbreak\ $\{0,1,34\}$,\allowbreak\ $\{0,1,36\}$,\allowbreak\ $\{0,1,39\}$,\allowbreak\ $\{0,1,40\}$,\allowbreak\ $\{0,1,41\}$,\allowbreak\ $\{0,1,46\}$,\allowbreak\ $\{0,1,49\}$,\allowbreak\ $\{0,1,\infty\}$,\allowbreak\ $\{0,2,10\}$,\allowbreak\ $\{0,2,12\}$,\allowbreak\ $\{0,2,15\}$,\allowbreak\ $\{0,2,23\}$,\allowbreak\ $\{0,2,25\}$,\allowbreak\ $\{0,2,28\}$,\allowbreak\ $\{0,2,42\}$,\allowbreak\ $\{0,2,44\}$,\allowbreak\ $\{0,5,19\}$,\allowbreak\ $\{0,5,26\}$,\allowbreak\ $\{0,5,34\}$,\allowbreak\ $\{0,5,39\}$,\allowbreak\ $\{0,5,51\}$,\allowbreak\ $\{0,5,\infty\}$.
\par}

\section{Classes of cross-ratio values}\label{app:tables}

For each colouring of Section~\ref{sec:line} we list the classes on which
$\tau$ takes the value $0$, each by its smallest element as an integer in
$\{2,\dots,q-1\}$. On every other class $\tau$ is $1$. For $F_{200}$ the
function is $\sigma$ and the classes are the narrow classes of
Section~\ref{sec:line}.

\paragraph{$D_{194}$: $q=193$, $15$ of $33$ classes.}
{\small\raggedright
$2$,\allowbreak\ $3$,\allowbreak\ $5$,\allowbreak\ $6$,\allowbreak\ $9$,\allowbreak\ $11$,\allowbreak\ $12$,\allowbreak\ $13$,\allowbreak\ $15$,\allowbreak\ $20$,\allowbreak\ $21$,\allowbreak\ $22$,\allowbreak\ $31$,\allowbreak\ $40$,\allowbreak\ $68$.
\par}

\paragraph{$F_{200}$: $q=199$, the group $\mathrm{PSL}(2,199)$, $20$ of $43$ narrow classes.}
{\small\raggedright
The function $\sigma$ takes the value $0$ on the narrow classes of the
following values $d$, so that colour $0$ is the union of the orbits of the
$4$-subsets $\{\infty,0,1,d\}$:
$d=2$,\allowbreak\ $3$,\allowbreak\ $7$,\allowbreak\ $11$,\allowbreak\ $12$,\allowbreak\ $13$,\allowbreak\ $15$,\allowbreak\ $17$,\allowbreak\ $24$,\allowbreak\ $27$,\allowbreak\ $29$,\allowbreak\ $35$,\allowbreak\ $37$,\allowbreak\ $38$,\allowbreak\ $41$,\allowbreak\ $55$,\allowbreak\ $68$,\allowbreak\ $73$,\allowbreak\ $81$,\allowbreak\ $107$.
\par}

\paragraph{$E_{500}$: $q=499$, $41$ of $84$ classes.}
{\small\raggedright
$3$,\allowbreak\ $5$,\allowbreak\ $6$,\allowbreak\ $9$,\allowbreak\ $10$,\allowbreak\ $11$,\allowbreak\ $12$,\allowbreak\ $14$,\allowbreak\ $15$,\allowbreak\ $17$,\allowbreak\ $19$,\allowbreak\ $20$,\allowbreak\ $24$,\allowbreak\ $27$,\allowbreak\ $32$,\allowbreak\ $34$,\allowbreak\ $39$,\allowbreak\ $43$,\allowbreak\ $47$,\allowbreak\ $56$,\allowbreak\ $63$,\allowbreak\ $67$,\allowbreak\ $71$,\allowbreak\ $72$,\allowbreak\ $75$,\allowbreak\ $77$,\allowbreak\ $82$,\allowbreak\ $92$,\allowbreak\ $110$,\allowbreak\ $115$,\allowbreak\ $116$,\allowbreak\ $117$,\allowbreak\ $119$,\allowbreak\ $124$,\allowbreak\ $144$,\allowbreak\ $168$,\allowbreak\ $172$,\allowbreak\ $177$,\allowbreak\ $191$,\allowbreak\ $223$,\allowbreak\ $224$.
\par}

\paragraph{$E_{600}$: $q=599$, $51$ of $100$ classes.}
{\small\raggedright
$3$,\allowbreak\ $4$,\allowbreak\ $5$,\allowbreak\ $10$,\allowbreak\ $12$,\allowbreak\ $14$,\allowbreak\ $15$,\allowbreak\ $16$,\allowbreak\ $18$,\allowbreak\ $20$,\allowbreak\ $22$,\allowbreak\ $23$,\allowbreak\ $24$,\allowbreak\ $33$,\allowbreak\ $34$,\allowbreak\ $35$,\allowbreak\ $37$,\allowbreak\ $38$,\allowbreak\ $39$,\allowbreak\ $42$,\allowbreak\ $43$,\allowbreak\ $49$,\allowbreak\ $53$,\allowbreak\ $55$,\allowbreak\ $65$,\allowbreak\ $66$,\allowbreak\ $67$,\allowbreak\ $69$,\allowbreak\ $72$,\allowbreak\ $79$,\allowbreak\ $81$,\allowbreak\ $84$,\allowbreak\ $85$,\allowbreak\ $90$,\allowbreak\ $95$,\allowbreak\ $103$,\allowbreak\ $104$,\allowbreak\ $106$,\allowbreak\ $112$,\allowbreak\ $115$,\allowbreak\ $127$,\allowbreak\ $138$,\allowbreak\ $139$,\allowbreak\ $143$,\allowbreak\ $146$,\allowbreak\ $147$,\allowbreak\ $174$,\allowbreak\ $197$,\allowbreak\ $199$,\allowbreak\ $241$,\allowbreak\ $258$.
\par}

\section{Where the searches found colourings}\label{app:census}

This appendix records the primes $q$ at which our searches found colourings
of $L_q$ that are invariant under the group named. Unless it is said
otherwise, the answers are those of a SAT solver, which decided each of
these small formulas in both directions. No proof certificate was recorded
for the answers UNSAT. They concern only the invariant family named and
establish no upper bound.

\paragraph{$(5,6;4)$, the group $\mathrm{PGL}(2,q)$.}
Among the $42$ primes $q$ with $101\le q\le 331$, the solver found a
colouring for $q=101$, $103$, $107$, $109$, $113$, $127$, $131$, $137$,
$139$, $151$, $157$, $163$, $179$, $191$ and $193$, and returned UNSAT for
the other $27$. For $q=193$ it reported that the colouring is unique.

\paragraph{$(5,6;4)$, the group $\mathrm{PSL}(2,q)$.}
Among the $30$ primes $q$ with $163\le q\le 331$, the solver found a
colouring for $q=163$, $167$, $179$, $191$, $193$ and $199$, and returned
UNSAT for the other $24$. For $q=167$ and $q=199$ the larger group
$\mathrm{PGL}(2,q)$ admits none.

\paragraph{$(6,6;4)$, the group $\mathrm{PGL}(2,q)$.}
Colourings were found for every prime $q$ with $7\le q\le 151$, and for
$q=157$, $199$, $251$, $317$, $331$, $337$, $347$, $349$, $359$, $397$,
$449$, $499$ and $599$; the one for $q=599$ was found by local search. For
$q=701$ neither the local search nor the solver gave an answer in the time
allowed.

\paragraph{$(5,5;4)$, the group $\mathrm{PGL}(2,q)$.}
Among the primes $q$ with $7\le q\le 199$, the solver found a colouring for
$q=7$, $13$, $17$, $23$ and $47$, each time unique up to exchanging the two
colours, and returned UNSAT for the others. The case $q=47$ is a colouring
of order $48$; it agrees with the bound $R(5,5;4)\ge 49$ of~\cite{ozdemir},
and we make no claim for it.
\phantomsection


\addcontentsline{toc}{section}{References}
\begin{thebibliography}{99}
\small

\bibitem{survey}
S.~P. Radziszowski,
\emph{Small Ramsey numbers},
Electron.\ J.\ Combin., Dynamic Survey DS1, revision \#18 (24 April 2026),
\url{https://doi.org/10.37236/21}.

\bibitem{mr1}
B.~D. McKay and S.~P. Radziszowski,
\emph{The first classical Ramsey number for hypergraphs is computed},
Proceedings of the Second Annual ACM--SIAM Symposium on Discrete Algorithms
(SODA '91), San Francisco, 1991, 304--308.

\bibitem{dyb2}
J.~Dybizba\'nski,
\emph{A lower bound on the hypergraph Ramsey number $R(4,5;3)$},
Contrib.\ Discrete Math.\ \textbf{13} (2018), no.~2, 112--115.

\bibitem{dyb3}
J.~Dybizba\'nski,
personal communication to S.~P. Radziszowski (2018), cited in~\cite{survey};
addendum to~\cite{dyb2} at \url{https://inf.ug.edu.pl/ramsey}.

\bibitem{brbh}
M.~Bruce, M.~Budden, and J.~Hiller,
\emph{Lexicographic products of $r$-uniform hypergraphs and some applications
to hypergraph Ramsey theory},
Australas.\ J.\ Combin.\ \textbf{70} (2018), 390--401.

\bibitem{yang}
H.~Yang,
\emph{Three improved lower bounds for hypergraph Ramsey numbers},
data deposit, Zenodo, 2026,
\url{https://doi.org/10.5281/zenodo.24301648}.

\bibitem{ozdemir}
K.~\"Ozdemir,
\emph{Lower bounds for the hypergraph Ramsey numbers $R(5,5;4)$ and
$R(5,6;4)$ from $\mathrm{PSL}(2,q)$-invariant colourings},
data deposit, Zenodo, 27 October 2026,
\url{https://doi.org/10.5281/zenodo.23689461}.

\bibitem{pysat}
A.~Ignatiev, A.~Morgado, and J.~Marques-Silva,
\emph{PySAT: A Python toolkit for prototyping with SAT oracles},
in: Theory and Applications of Satisfiability Testing -- SAT 2018,
Lecture Notes in Comput.\ Sci.\ 10929, Springer, 2018, 428--437.

\bibitem{glucose}
G.~Audemard and L.~Simon,
\emph{Predicting learnt clauses quality in modern SAT solvers},
in: Proceedings of the 21st International Joint Conference on Artificial
Intelligence (IJCAI 2009), 399--404.

\bibitem{cadical}
A.~Biere, K.~Fazekas, M.~Fleury, and M.~Heisinger,
\emph{CaDiCaL, Kissat, Paracooba, Plingeling and Treengeling entering the SAT
Competition 2020},
in: Proceedings of SAT Competition 2020 -- Solver and Benchmark Descriptions,
Department of Computer Science Report Series B, vol.~B-2020-1, University of
Helsinki, 2020, 51--53.

\bibitem{sonyl}
E.~Song, W.~Ye, and Y.~Liu,
\emph{New lower bounds for Ramsey number $R(p,q;4)$},
Discrete Math.\ \textbf{145} (1995), 343--346.

\bibitem{liu}
S.~C. Liu,
\emph{Lower bounds for small Ramsey numbers on hypergraphs},
preprint, 2019, arXiv:1906.00132.

\bibitem{probsat}
A.~Balint and U.~Sch\"oning,
\emph{Choosing probability distributions for stochastic local search and the
role of make versus break},
in: Theory and Applications of Satisfiability Testing -- SAT 2012,
Lecture Notes in Comput.\ Sci.\ 7317, Springer, 2012, 16--29.

\end{thebibliography}
\end{document}